\documentclass{article}
\usepackage{amsthm,amsmath,amssymb,enumerate}
\usepackage[left=1in,top=1.2in,right=1in,bottom=1.2in,head=.5in,foot=.3in]{geometry}

\usepackage{fancyhdr}

\newtheorem{thm}{Theorem}

\theoremstyle{definition}
\newtheorem{defn}{Definition}
\newtheorem{example}{Example}
\newtheorem{rem}{Remark}

\usepackage[]{graphicx}

\DeclareMathOperator{\Div}{Div}
\newcommand{\lx}{^{\leq x}}

\newcommand{\A}{\mathcal A}

\newcommand{\M}{\mathcal M}
\newcommand{\N}{\mathcal N}
\renewcommand{\O}{\mathcal O}
\newcommand{\U}{\mathcal U}
\newcommand{\bbN}{\mathbb N}
\newcommand{\K}{\mathbb K}
\newcommand{\Q}{\mathbb Q}
\newcommand{\R}{\mathbb R}
\newcommand{\Z}{\mathbb Z}
\newcommand{\Ltf}{\L_{\text{tf}}}
\DeclareMathOperator{\prov}{prov}
\DeclareMathOperator{\Con}{Con}
\DeclareMathOperator{\Th}{Th}

\let\losell\L
\newcommand{\los}{\losell os}
\renewcommand{\L}{\mathcal L}

\begin{document}

\title{Computable Axiomatizability of Elementary Classes}
\author{Peter Sinclair}
\date{}

\maketitle

\begin{abstract}
  The goal of this paper is to generalise Alex Rennet's proof of the non-axiomatizability of the class of pseudo-o-minimal structures. Rennet showed that if $\L$ is an expansion of the language of ordered fields and $\K$ is the class of pseudo-o-minimal $\L$-structures ($\L$-structures elementarily equivalent to an ultraproduct of o-minimal structures) then $\K$ is not computably axiomatizable. We give a general version of this theorem, and apply it to several classes of topological structures.
\end{abstract}

\section{Introduction}

Given a class $\K$ of $\L$-structures, we write $\Th(\K)$ for the first order theory of $\K$; that is, the set of all $\L$-sentences that are true in every structure of $\K$. Recall that a class $\K$ is called \emph{elementary} when $\M\models\Th(\K)$ if and only if $\M$ is an element of $\K$, and that this holds if and only if $\K$ is closed under ultraproducts and ultraroots \cite[Corollary 8.5.13]{Hod}. We say that an elementary class $\K$ is \emph{computably axiomatizable} if there is a computable axiomatization of $\Th(\K)$. With this terminology, Rennet proved that the class of pseudo-o-minimal fields (fields which are elementarily equivalent to an ultraproduct of o-minimal structures) is not computably axiomatizable \cite{Ren-P}.

Rennet's paper was motivated by a number of results, among them Ax's proof \cite{Ax} that the theory of finite fields is decidable, and hence that the class of pseudo-finite fields is computably axiomatizable. As with the class of finite fields in the language of rings, the class of o-minimal structures in a language with an ordering and an extra unary predicate is not elementary. For each $n\in\bbN$, let $\M_n$ be a copy of the real numbers in this language, where the ordering is interpreted by the usual ordering and the unary predicate is interpreted as $\{0,1,\ldots,n\}$. It is easy to see that each $\M_n$ is o-minimal, but that the ultraproduct has a copy of the natural numbers as a definable set; this is clearly not a finite union of points and intervals, and hence the ultraproduct is not o-minimal. Thus, the class of o-minimal structures is not closed under ultraproducts, and so is not elementary.

Multiple proposals were made for possible axiomatizations of the class of pseudo-o-minimal structures (see \cite{Fo} and \cite{Sch}, for instance). However, Rennet showed that in the case where the language expands that of ordered fields, there is no computable axiomatization for the theory of o-minimality, and hence the class of pseudo-o-minimal structures is not computably axiomatizable.

In \cite{HM-C}, Haskell and Macpherson developed the notion of $C$-minimality, a generalization of o-minimality obtained by replacing the binary ordering by a ternary relation. Haskell and Macpherson looked at another generalization of o-minimality in \cite{HM-P}, $P$-minimality, which is defined so that $P$-minimal fields are $p$-adically closed, just as o-minimal fields are real closed. Given the similarities between these settings and o-minimality, they are both contexts in which it is natural to ask whether Rennet's theorem applies.

In this paper, we adapt Rennet's proof to give a more general theorem, which can then be applied to other classes, including those of $C$-minimal and $P$-minimal structures. Section 2 contains the preliminaries and proof of the generalized theorem, while Section 3 contains some examples, including those mentioned above.

\section{Preliminaries and the Generalized Theorem}

We state our generalization of Rennet's theorem in the context of first order topological stuctures, as introduced by Pillay in \cite{Pil}:

\begin{defn} Let $\A$ be a structure in a language with a formula $B(x,\bar y)$ (where $x$ is a single variable and $\bar y$ is a tuple) such that the set of $A$-subsets $\{B(x,\bar a)^\A :\ \bar a\subseteq A\}$ is a basis for a topology on $A$. We say that such an $\A$ is a \emph{first-order topological structure}, or simply a topological structure. Note that for any $\A'\equiv\A$, $(\A',B)$ is also a topological structure.

We extend this notion by saying a class $\K$ of $\L$-structures is \emph{uniformly topological} if there is a single formula $B$ such that each $\A\in\K$ is a topological structure with a basis given by $B$.
\end{defn}

Recall the notion of a provability relation which plays a fundamental role in the proof of G\"odel's Second Incompleteness Theorem (see, for instance, \cite{BBJ}): if $\Gamma$ is a computable list of sentences in the language of arithmetic then there exists a binary relation $\prov(s,d)$ such that in the standard model of Peano Arithmetic, $\prov(s,d)$ if and only if $d$ is the code number of a sentence and $s$ is the code number for a proof of that sentence from $\Gamma$.

\begin{thm}\label{main} Fix any computable language $\L$ containing a unary predicate $N$. Suppose $\K$ is a uniformly topological class of $\L$-structures whose topology is given by the formula $B(x,\bar y)$. Moreover, suppose that for each $\A\in\K$, discrete definable subsets of $A$ are finite. Let $\Lambda$ be any computable subset of $\Th(\K)$.

Fix distinguished $\L$-formulas $\alpha$, $\mu$, and $\leq$ which define subsets of $N^3$, $N^3$, and $N^2$, respectively, without parameters. Also fix $\emptyset$-definable constants $0,1\in N$. Let $T$ be the $\L$-theory described below:
\begin{enumerate}[(I)]
\item $(N,\alpha,\mu,\leq,0,1)$ is a model of the relational theory of Peano Arithmetic, $PA$.
\item $N$ is discrete: that is, $T$ contains the sentence
\[ \forall x\in N\ \exists \bar a \forall y (y\in N\wedge B(y, \bar a) \rightarrow y = x) .\]
\item For each $\psi\in\Lambda$, $T$ contains $\forall x\in N\ \psi\lx$, where $\psi\lx$ is the sentence $\psi$ with any occurrence of $N(t)$ replaced by $N(t) \wedge t \leq x$.
\end{enumerate}
If $T$ is consistent then there is an $\L$-structure $\mathcal R^{\L}_\Lambda$ which satisfies $\Lambda$, but is not elementarily equivalent to an ultraproduct of structures in $\K$.
\end{thm}

It follows that the class $\{\M:\ \M\models \Th(\K)\}$ is not computably axiomatizable, since given any potential axiomatization $\Lambda$, the structure $R^{\L}_\Lambda$ obtained in the theorem satisfies $\Lambda$ but not $\Th(\K)$.

\begin{proof} Assume that $T$ is consistent. In every model of $T$, the interpretation of $N$ is a model of Peano Arithmetic, and so by G\"odel's Second Incompleteness Theorem, $T+\neg\Con(T)$ is also consistent. Thus, there exists a model $\A$ of $T+\neg\Con(T)$. In particular, if $\prov(s,d)$ is the provability relation for $T$ and $c$ is the G\"odel number for the statement $0=1$ then $\A \models \exists s \prov(s,c)$; that is, there exists $a\in N$ with $\A \models \prov(a,c)$.

Fix $x\in N$ with $x$ sufficiently large to code the proof of $c$ (among other conditions, $a\leq x$ and $c\leq x$) and consider the structure $\A_x$ which is identical to $\A$ except that $N$ is replaced by the initial segment $\{n\in N_\A :\ n\leq x\}$. Since $\A$ satisfies the axiom schema (III), $\A_x$ satisfies $\Lambda$. By Theorem 2.7 of \cite{Ka91}, since $N_{\A_x}$ is an initial segment of $N_\A$, a model of the relational theory of Peano Arithmetic, it is a $\Delta_0$-elementary substructure of $N_\A$. Thus, since $a$ being a code for a proof of $0=1$ in $T$ is a $\Delta_0$-property of $a\in N_{\A_x}$, we have $N_{\A_x} \models \exists s \prov(s,c)$.

We claim that $\A_x$ is the desired structure $\mathcal R^{\L}_\Lambda$. Suppose for contradiction that $\A_x$ is elementarily equivalent to an ultraproduct of structures in $\K$:
\[ \A_x \equiv \A' = \prod_{i\in I} \A_i/\U \]
where $\U$ is a non-principal ultrafilter on $I$, and every $\A_i$ is a structure in $\K$. Since property (II), that $N$ is discrete, is described by a first order sentence, it also holds in $\A'$, and hence, by \los's Theorem, it also holds in $\U$-most of the $\A_i$. Since each $\A_i\in\K$ and each $N_{\A_i}$ is trivially definable, by assumption $\U$-most of the $N_{\A_i}$ are finite.

Then, since $\N_{\A_x}$ is an initial segment of a model of $PA$, so is $N_{\A'}$ and $\U$-most of the $N_{\A_i}$. But $\U$-most of the $N_{\A_i}$ are finite, so $\U$-most of the $N_{\A_i}$ are finite initial segments of a model of $PA$, and hence are isomorphic to a substructure of $\bbN$ with universe $I_n = \{0,1,\ldots,n\}$ for some $n\in\bbN$. That is, $\U$-most $N_{\A_i}$ are isomorphic, for some $n_i$, to the structure
\[ \bbN_{n_i} = (I_{n_i},\ \{(x,y,z)\in I_{n_i}: x+y=z\},\ \{(x,y,z)\in I_{n_i}: xy=z\},\ \{(x,y)\in I_{n_i}: x\leq y\}) . \]

Let $c'\in N_{\A'}$ be a code for $0=1$ and $\prov(d,s)$ the provability relation for $T$. Since $N_{\A_x} \equiv N_{\A'}$, we have $N_{\A'}\models\exists s \prov(s,c')$. Choose an index $i$ such that $N_{\A_i}\models \exists s\prov(s,c_i')$ and $N_{\A_i}$ is isomorphic to some $\bbN_{n_i}$ as above. Then, since $N_{\A_i}\cong \bbN_{n_i}$ is a $\Delta_0$-elementary substructure of $\bbN$, there exists $b\in\bbN$ such that $\bbN\models \prov(b,c)$, where $c\in\bbN$ is the image of $c_i'\in N_{\A_i}$. Because of the interpretation of $\prov(b,c)$ in the standard model $\bbN$, this $b$ corresponds to an actual proof of $0=1$ in $T$. Hence $T$ is inconsistent, contradicting our assumption, and so $\A_x$ cannot be elementarily equivalent to an ultraproduct of structures in $\K$.
\end{proof}

\begin{rem} Note that the requirement of the predicate $N$ being included in the language is merely a convenience. Any occurrence of $N$ could be replaced by a distinguished formula in one variable and the proof would be unaffected.
\end{rem}

\section{Consequences}

The examples below are all straightforward consequences of the theorem, which amount to choosing an appropriate class for $\K$ and showing that the theory $T$ from the theorem is consistent.

The first pair of examples, $P$-minimality and $C$-minimality, are variations of o-minimality designed for valued fields. While more detailed descriptions can be found in \cite{HM-P} and \cite{HM-C}, for our purposes we need only a single example of each to use in our construction of a model of $T$.

Fix a prime $p$. Then any rational number can be written in the form $p^n\frac{a}{b}$ where $n,a,b\in\Z$ and $p\nmid a,b$. We define a valuation $v_p: \Q\to\Z$ by $v_p(p^n\frac{a}{b}) = n$. With appropriate choices of language, the completion $\Q$ with respect to the norm $|x| = p^{-v(x)}$ is an example of a $P$-minimal structure, denoted $\Q_p$. In Chapter III of \cite{Kob}, Koblitz shows that the metric completion of the algebraic closure of $\Q_p$, denoted $\Omega_p$, is an algebraically closed valued field; it then follows from \cite[Theorem C]{HM-C} that $\Omega_p$ is an example of a $C$-minimal structure.

Let $K$ be one of the fields described in the previous paragraph. In both cases, the exponential function $\exp(x) = \sum_{n=0}^\infty \frac{x^n}{n!}$ converges on the set $p\O = \{px\in K:\ v(x)\geq 0\}$, and is bijective on this domain. Moreover, $\Q_p$ and $\Omega_p$ continue to be examples of $P$-minimal and $C$-minimal structures when the language is expanded by adding a symbol for the exponential function restricted to $p\O$; see \cite[Theorem B]{DHM} for the $P$-minimal case and \cite[Theorem 1.6]{LR} for the $C$-minimal case.

We create a model of Peano arithmetic in $K$ as follows: take $N = \{p^{pn}:\ n\in\bbN\}$, and define $\{0_N,1_N,\alpha,\mu,\leq\}$ via the natural bijection $p^{pn} \mapsto n$. Note that these sets will not be definable in $K$ using the usual language for $P$-minimal or $C$-minimal fields, even after adding a symbol for the restricted exponential function. Clearly, $N$ will be isomorphic to the usual interpretaion of the natural numbers, and hence will be a model of Peano arithmetic. In the examples below, we simply need to show that this structure is definable in our chosen language; the additional factor of $p$ in the exponent will be required to ensure that $\exp(x)$ is defined everywhere required.

\begin{example} Let $\L_d = \{+,-,\cdot,0,1,\operatorname{Div},\{P_n\}_{n\in\N}\}$ be the language used in \cite{HM-P}, let $\L$ be any expansion of $\L_d\cup\{\exp,N\}$, and let $\K$ be the class of $P$-minimal $\L$-structures in which $\exp$ is interpreted as the restricted exponential. Then the class $\K' = \{\A:\ \A\models\Th(\K)\}$ is not computably axiomatizable.
\end{example}

\begin{proof} Let $\Lambda$ be a purported axiomatization of $\Th(\K)$, and note that each $\A\in\K$ has a topology with uniformly definable basis $B(x,c,d) = \{x\in A :\ \Div(x-c,d) \wedge \neg(x-c=d)\}$. It follows from Lemma 4.3 of \cite{HM-P} that every discrete definable set in a $P$-minimal structure is finite.

To show $T$ is consistent, we consider $\Q_p$ with $N = \{p^{pn}:\ n\in\bbN\}$ and $\{0_N,1_N,\alpha,\mu,\leq\}$ interpreted as described above. Clearly, $0_N$ and $1_N$ are $\emptyset$-definable, and $x\leq y$ is equivalent to $\Div(x,y)$. Moreover, $\alpha(x,y,z)$ is defined by $xy=z$. It remains to show that $\mu(x,y,z)$ is definable in the language.

As noted above, the restricted exponential function on $\Q_p$ is bijective, and hence the function $\ln(x)$ given by $\ln(x) = y$ when $\exp(y) = x$ is definable in $\L$ for $v(x)\geq 1$. We can thus take $\mu(x,y,z)$ to be the set defined by
\[ \exp\left(\frac{\ln(x)\ln(y)}{p^2\ln(p)} \right) = z . \]

To turn this into an $\L$-structure $\A$, we simply use a trivial interpretation of every relation, function, and constant symbol not in $\L_d\cup\{\exp,N\}$. Conditions (I) and (II) for $T$ are satisfied by choice of $N$. Condition (III) follows from the fact that every initial segment of $N$ is finite: for all $x\in N$, $\A_x$ is a definitional expansion of $\Q_p$ (as an $\L_d$-structure), which means it is $P$-minimal, and hence satisfies $\Lambda$. Thus $\A\models T$, and so by Theorem \ref{main}, there is a model of $\Lambda$ which is not an element of $\K'$.
\end{proof}

\begin{example} Let $\L_c = \{+,-,\cdot,0,1,C\}$ be the language of $C$-minimal fields described in \cite{HM-C}, let $\L$ be any proper expansion of $\L_c\cup\{\exp\}$, and let $\K$ be the class of $C$-minimal $\L$-structures in which $\exp$ is interpreted as the restricted exponential. Then the class $\K' = \{\A:\ \A\models\Th(\K)\}$ is not computably axiomatizable.
\end{example}

\begin{proof} Let $\Lambda$ be a purported axiomatization of $\Th(\K)$, and note that $B(x,b,c) = \{x:\ C(b;x,c)\}$ gives a uniformly definable basis for a topology on each $\A\in\K$. As noted in Lemma 2.4 of \cite{HM-C}, discrete definable sets in $C$-minimal structures are finite.

To show $T$ is consistent, consider $\Omega_p$ with $N = \{p^{pn}:\ n\in\bbN\}$ and $\{0_N,1_N,\alpha,\mu,\leq\}$ interpreted as described above. Again, $0_N$ and $1_N$ are $\emptyset$-definable, and $x\leq y$ is equivalent to $\neg C(y;x,1)$, where $1$ here is the multiplicative identity in the field, not in the set $N$. The exponential function is again bijective, which means $\alpha$ and $\mu$ are definable by the same formulas as in the $P$-minimal case. Then we can form an $\L$-structure in the same way as before, and it will satisfy $T$ for the same reasons described above.
\end{proof}

For our final two examples, we look to Pillay's paper \cite{Pil}. In section 3 of that paper, Pillay defines a dimension rank $D_A$ for first order topological structures, which we will not repeat here. He notes that every stable first order topological structure has the discrete topology, and so Theorem \ref{main} cannot be applied to stable structures. However, he introduces a different notion of stability for such structures, which can be used:

\begin{defn} A first order topological structure $\A$ is said to be \emph{topologically totally transcendental}, or t.t.t., if it satisfies the following properties:
\begin{enumerate}[(A)]
\item Every definable set $X\subseteq A$ is a boolean combination of definable open sets.
\item Every definable set $X\subseteq A$ has $d(X) < \infty$, where $d(X)$ is the maximum choice of $d$ such that $X$ can be written as a disjoint union of nonempty definable sets $X_1,\ldots,X_d$ with each $X_i$ both closed and open in $X$.
\item $A$ has dimension, meaning $D_A(A) < \infty$.
\item The topology on $A$ is Hausdorff.
\end{enumerate}
Moreover, $\A$ is said to be \emph{$t$-minimal} if $\A$ is t.t.t.\ and $D_A(A) = d(A) = 1$. 
\end{defn}

In the case of an ordered structure, $t$-minimality is equivalent to o-minimality \cite[Proposition 6.2]{Pil}. However, the defition is less restrictive in general. Since the ordering on the reals is definable in the field language, $(\R,+,\cdot)$ with the usual topology is $t$-minimal, while the structure $(\mathbb C,+,\cdot,P)$ with the usual topology and $P$ interpreted as a predicate for the positive reals is an example of a t.t.t.\ structure which is not $t$-minimal.

\begin{example} Let $\L$ be a proper expansion of $\Ltf = \{+,\cdot,0,1,B\}$, where $+$ and $\cdot$ are binary function symbols and $B$ is an $n$-ary relation symbol for some $n\geq 2$, and let $\K$ be the class of $t$-minimal $\L$-structures in which $B(x,\bar y)$ gives a basis for a topology. Then the class $\K' = \{\A:\ \A\models\Th(\K)\}$ is not computably axiomatizable.
\end{example}

\begin{proof} Let $\Lambda$ be a purported axiomatization of $\Th(\K)$, and suppose $\A\in\K$ with $N\subseteq \A$ discrete and definable. Since $\A$ is Hausdorff, each point $a\in N$ is closed, and since $N$ is discrete, each $a\in N$ is open in $N$. Thus, $|N| = d(N)$ is finite by condition (B), and so discrete definable subsets in each $\A\in\K$ are finite.

Consider the real numbers with the usual interpretation of $+$, $\cdot$, 0, and 1, and $N$ as a predicate for the natural numbers. If $I$ is the set of all open intervals with endpoints in $\R$, then $|I| = |\R|$, so there exists a bijection $f: \R\to I$; take $B(x_1,\ldots,x_n)$ to be the relation $x_1\in f(x_2)$. Taking a trivial interpretation of every function, relation, and constant symbol not in $\Ltf$ gives an $\L$-structure $\A$, which we claim is a model of $T$.

For (I), take $0_N = 0$, $1_N = 1$, $\alpha$ and $\mu$ the graphs of $+$ and $\cdot$ restricted to $N$, and $x\leq y$ iff $x,y\in N$ and $\exists z (x+z^2 = y)$. Clearly, this gives a model of Peano Arithmetic. Since $N\cap (a-1,a+1) = \{a\}$ for every $a\in N$, we have (II), that $N$ is discrete. It remains to show that for any $x\in N$, the structure $\A_x$ is t.t.t.

First, note that $B$ gives the usual topology on $\R$, which is clearly Hausdorff, and thus we have condition (D) of t.t.t. Moreover, the definable sets in $\A_x$ are precisely the same as those in $(\R,+,\cdot,0,1,\leq)$, and hence are finite unions of points and intervals: this gives conditions (A) and (B). Finally, any definable set $X\subseteq A$ without interior in $A$ must be a finite union of points, in which case $D_A(X) = 0$, and so $D_A(A) = 1$. This is equivalent to condition (C) by Proposition 3.7 of \cite{Pil}. Thus, $\A$ satisfies condition (III), which means $T$ is consistent and Theorem \ref{main} can be applied.
\end{proof}

\begin{rem} As with $N$, the inclusion of $B$ in the language is merely a convenience. Given a distinguished formula for $B$ that satisfies the assumptions for the structure to be t.t.t., we could (with more difficulty) interpret the function and relation symbols in such a way that we obtain essentially the same model of $T$ given above.
\end{rem}

\begin{example} Let $\L$ be an expansion of $\{+,\cdot,0,1,B,N\}$, where $+$ and $\cdot$ are binary function symbols and $B$ is an $n$-ary relation symbol for some $n\geq 2$, and let $\K$ be the class of $t$-minimal $\L$-structures in which $B(x,\bar y)$ gives a basis for a topology. Then the class $\K' = \{\A:\ \A\models \Th(K)\}$ is not computably axiomatizable.
\end{example}

\begin{proof} In the previous example, we have already shown everything necessary except that the structure $\A_x$ has $d(A) = 1$. But this is equivalent to saying that $\R$ (with its usual topology) is connected, which is clearly true.
\end{proof}

\bibliographystyle{mlq}

\providecommand{\WileyBibTextsc}{}
\let\textsc\WileyBibTextsc
\providecommand{\othercit}{}
\providecommand{\jr}[1]{#1}
\providecommand{\etal}{~et~al.}

\end{document}